\documentclass{article}
\usepackage{latexsym}
\usepackage{amsmath}
\usepackage{graphicx}
\usepackage{amsfonts}
\usepackage{amssymb}
\usepackage{amsthm}

\newtheorem{theorem}{Theorem}

\newtheorem{algorithm}[theorem]{Algorithm}

\newtheorem{conjecture}[theorem]{Conjecture}
\newtheorem{corollary}[theorem]{Corollary}

\newtheorem{definition}[theorem]{Definition}

\newtheorem{lemma}[theorem]{Lemma}

\newtheorem{problem}[theorem]{Problem}
\newtheorem{proposition}[theorem]{Proposition}
\newtheorem{remark}[theorem]{Remark}

\DeclareMathOperator{\Char}{char} \DeclareMathOperator{\Hom}{Hom}
\DeclareMathOperator{\trunc}{trunc} \DeclareMathOperator{\Ad}{Ad}

\newcommand{\FM}{{\mathfrak m}}

\begin{document}
\author{J. Kuttler\footnote{Research partially supported by
    the Swiss National Science Foundation}~and N. R. Wallach\footnote{Research partially supported by a NSF Summer grant}}
\title{Representations of $SL_2$ and the distribution of points in
$\mathbb P^n$}
\date{November 13, 2002}

\maketitle

\section{Introduction} \label{intro}

In 1958 Nagata \cite{NAG59} gave an ingenious argument that
demonstrated the existence of counterexamples to Hilbert's
Fourteenth Problem. \ Recall that the original problem is the
following: Let $K = F(x_1,x_2,\dots,x_n)$ be the function field
of affine $n$-space $V = F^n$ over an algebraically
closed field $F$, and suppose $L \subset K$ is any subfield. 
Then the question is: Is $A = L \cap F[x_1,x_2,\dots,x_n]$ a
finitely generated $F$-algebra? In most cases of interest, $L$
is the field of invariants of an algebraic group $G$ acting
linearly on $V$, and $A$ becomes the ring of invariant regular
functions on $V$.  Certainly, if $G$ is reductive, the answer is
yes, a result due to Hilbert himself for $\Char(F) = 0$, but it is much
more subtle in positive characteristic. More
generally for $G$ reductive $\mathcal O(X)^G$ is a finitely generated
$F$-algebra for any affine variety $X$ over $F$ on which $G$
acts.  On the other hand, if $G = F^+$ is the additive group,
due to more recent work of for example Daigle and Freudenburg \cite{DF99}, there are
actions of $G$ on $V$ for certain values of $n = \dim V$ (e.g.\ $n =
5$), such that $\mathcal O(V)^G$ is not finitely generated, even in
characteristic zero.  By a well known result of Weitzenb\"ock 
these actions cannot be linear, since in characteristic zero linear actions of
$F^+$ always
extend to representations of $SL_2(F)$.  On the other hand it
was not until Nagata found his counterexample, that the general
question of finite generation of invariants for a representation
of an algebraic group was settled. The purpose of this paper is
to investigate a byproduct of the remarkably involved
methods which Nagata used to prove the following result.

\begin{theorem}[Nagata's counterexample] \label{counterex} There
 is a linear action of $(F^+)^{13}$ on $V = F^{32}$, such that the ring of invariant functions is not finitely
generated.
\end{theorem}

For quite some time $13$ and $32$ were the lowest known
dimensions. Recently, Steinberg \cite{STEIN} improved this
result to $6$ and $18$, i.e. $(F^+)^{6}$ acting on $F^{18}$, using similar
methods. Much of the work  in the present paper has been inspired by
Steinberg's result.

The idea of the proof of theorem~\ref{counterex} is essentially
the following: Nagata defined a sequence of ideals $\{\mathfrak
a_m\}_{m \geq 0}$ with $\mathfrak a_m\mathfrak a_n \subset
\mathfrak a_{m+n}$ of the polynomial ring in three variables
$F[x,y,z]$ which satisfies the following condition:  for each
positive integer $m$, there are other positive integers $n > m$
and $k$ such that $\mathfrak a^k_n \not= \mathfrak a_{kn}$. \
Moreover the invariants of the product of $G = (F^+)^{13}$ with a
certain torus $T = (F^*)^{16}$, centralizing $G$, are isomorphic
with $R = \bigoplus_{m \geq 0} \mathfrak a_mt^{-m} \subset
F[x,y,z,t,t^{-1}]$. Since $T$ is reductive, the invariants of $G$
are finitely generated if and only if the invariants of $G \times
T$ are. The cited condition above implies the failure of an Artin-Rees type lemma for $R$, hence $R$ is not finitely generated.
 The ideal $\mathfrak a_m$ is the set of polynomials in
$F[x,y,z]$, vanishing at $16$ generic lines in $F^3$ with
multiplicity $m$. As one sees from this discussion, certain
properties of the ideals $\mathfrak a_m$ are essential to
Nagata's argument. He himself devoted considerable effort to the
general problem of this paper in dimension $2$: to determine the (dimension of the)
ideal of homogeneous polynomials vanishing at finitely many
points in $\mathbb P^2$ to a certain multiplicity, and
algorithmically solved the problem for less than nine points. We will present below various aspects of this problem, including
our version of his algorithm.

The first named author would like to thank the University of California San
Diego for its hospitality during his work on this paper. 

\section{The general problem} \label{general}
Let $F$ be an algebraically closed field. Let $V = F^{n+1}$ be affine $(n+1)$-space over $F$, and denote the
projective space of $V$ by $\mathbb P(V)$. Let $\mathcal O(V)$ be
the space of regular functions (i.e.\ polynomials) on $V$,
equipped with the usual grading $\mathcal O(V) = \bigoplus_{d\geq
0} \mathcal O(V)_d$. Similarly, $\mathcal D(V)$ is the graded
algebra of differential operators on $\mathcal O(V)$ with values
in $\mathcal O(V)$ and constant coefficients.

Now let $p \in \mathbb P(V)$ be arbitrary. For a positive integer
$m$, we say a homogeneous function $f \in \mathcal O(V)$ vanishes
to order $m$ at $p$, if $f \in \FM_p^m$, where $\mathfrak{m}_p
\subset \mathcal O(V)$ is the homogeneous prime ideal of height
one defining $p$.

\begin{remark} \label{vanishingorder} Throughout the paper, we will use several different
notions of vanishing to a given order at a point $p$, which appear
in the following list. Suppose $f \in \mathcal O(V)$ is a
homogeneous polynomial. The following statements are equivalent:
\medskip
\begin{itemize}
\item[i)] $f$ vanishes to order $m$ at $p$.
\item[ii)] $f$ vanishes to order $m$ at one point $\tilde p \in p
\setminus\{0\} \subset V$
\item[iii)] $f$ vanishes to order $m$ at all points of $p \setminus\{0\} \subset
V$
\item[iv)] $f$ vanishes to order $m$ at all points of $p \subset
V$.
\item[v)] The hypersurface $H \subset \mathbb P(V)$ defined by $f$
contains $p$ with multiplicity $m$.
\item[vi)] If $\Char(F) = 0$, then for every differential operator $D \in \mathcal D(V)$ of
order less or equal to $m-1$, we have $D(f) \in \FM_p$. Or, using
$ii)$ again, for every such $D$, $D(f)$ vanishes at a fixed
$\tilde p \in p\setminus \{0\}$.
\end{itemize}
\end{remark}

We will now describe the general problem: For $l$ points
$p_1,p_2,\dots,p_l \in \mathbb P(V)$ and nonnegative integers
$m_1,m_2,\dots,m_l$ and $k$ we set $\mathbf p = (p_1,\dots,p_l)$
and $\mathbf m = (m_1,\dots,m_l)$. Let
$$d(\mathbf p;\mathbf m;k) = \dim_F \bigcap_i
\FM_{p_i,k}^{m_i}$$ where $\FM_{p_i,k}^{m_i}$ denotes the degree
$k$-part of $\FM_{p_i}^{m_i}$. In other words, $d(\dots)$ denotes
the dimension of the space of all homogeneous polynomials of
degree $k$ that vanish at each of the $p_i$ to order $m_i$. Let
$d(\mathbf m;k)$ denote the minimum over $d(\mathbf p; \mathbf
m;k)$ where the $p_1,p_2,\dots,p_l$ vary. $l$ points will be called
\emph{generic} if $d(\mathbf p; \mathbf m; k) = d(\mathbf m;k)$ for all
$\mathbf m$. We will be studying
these numbers. The ideal of functions
vanishing of order $m_i$ at $p_i$ is denoted $I_{\mathbf
p,\mathbf m}$ or also $I_{\mathbf m}$, with the points
understood. The problem of determining the numbers $d(\mathbf
m;k)$ is well known, and a substantial portion of the known results are due to
Nagata, as we mentioned in the introduction. In his famous
papers on rational surfaces \cite{NAG60_I},\cite{NAG60_II} he
developed among other things a theory on linear systems,
generalizing the ordinary terminology of linear systems, which
allowed him to give a complete algorithm to determine these
numbers for $n = 2$ and $l \leq 9$. We will give an alternate
proof of this algorithm using representation theory of $SL_2$. It
is remarkable that up today, no general solution to the problem
seems to be known known except for some special cases.

There are several questions related to our problem, whose answer
would lead to a partial solution. We will now assume that
$\Char(F) = 0$.

\begin{problem} \label{independence}
When does `vanishing at order $m$ at $l$ points' impose
independent conditions on $\mathcal O(V)_k$?
\end{problem}

Of course, independence here refers to the question, whether the
obvious equations are independent: for each point $p$, vanishing
to order $m$ at $p$ is the same as solving the equations $D(f)(p)
= 0$ for all differential operators $D \in \mathcal D(V)$ of
degree less or equal $m-1$. These equations certainly never are
independent. But if we choose a complement to the trivial
equations, namely the elements of the ideal $D_p\mathcal D(V)$,
where $D_p$ means derivation along $p$, then the question makes
sense. Note that this ideal is independent of the choice of a
representative of $p$. So the number of equations is $\sum_{i=1}^{l}
\binom{m_i - 1 + {n}}{n}$, and the problem asks, when is this the
codimension of $I_{\mathbf m;k}$ in $ \mathcal O(V)_k$. we also
say that the conditions are independent, if they cut out
everything, i.e.\ if $I_{\mathbf m;k} = 0$. We will come back to
this later, but note for now, it is easy to see, that for $k$
large enough the answer to problem~\ref{independence} is yes (in
any characteristic).

From a more algebro-geometric point of view, one can think of our
$l$ generic points with multiplicities $m_1,\dots,m_l$ as a zero
dimensional subscheme $Z$ of $\mathbb P(V)$, defined by the sheaf
of ideals $\mathcal I  = \bigcap \mathcal I_{p_i}^{m_i}$, where
$\mathcal I_{p_i}$ is the ideal sheaf defining the reduced point
$p_i$. Then one has the natural exact sequence
\begin{equation}
0 \rightarrow \mathcal I \rightarrow \mathcal O_{\mathbb P(V)}
\rightarrow \mathcal O_Z \rightarrow 0
\end{equation}
Since $Z$ is a zero-dimensional scheme, we may write $\mathcal
O_Z = \bigoplus_i \mathcal O_{Z,p_i} = \bigoplus_i  F^{b_i}$,
with $b_i = \binom{m_i - 1 + n}{n}$. And our problem~\ref{independence} becomes equivalent to the following:

\begin{problem}\label{independence2} For which $k$ has the induced map
\begin{equation}
H^0(\mathbb P(V),\mathcal O_{\mathbb P(V)}(k)) \rightarrow
H^0(Z,\mathcal O_Z(k))
\end{equation}
maximal rank?
\end{problem}

Since we know, that for large $k$ $H^1(\mathbb P(V),\mathcal I(k))
= 0$, our claim above follows: the answer is yes, for $k$
sufficiently large.

In this form, the problem has been studied by Alexander and
Hirschowitz, for example. Alexander showed in \cite{A} that
for $k \geq 5$ and all $m_i$ equal to $2$, the codimension of
$I(k)$ is either $l(n+1)$, or $I(k) = 0$. In other words the rank
in problem~\ref{independence2} is always maximal.

Alexander and Hirschowitz \cite{AH92} also study the
case $k = 4$. For $k = 4$ and  all $m_i = 2$
the same is true, except for the following cases:
$$n=2,l=2; \qquad n=3,l=9; \qquad n=4, l=14.$$

Moreover, recently Ciliberto and Miranda \cite{CM2000} showed that whenever $n
= 2$ and all multiplicities equal and less than or equal $12$, then problem~\ref{independence} has an affirmative answer for $k \geq 3m$.

To conclude this section we return to the original question of
Nagata. In his definition of $\mathfrak a_m$ all the $m_i$s are
equal to $m$. So as a part of our problem, one might ask

\begin{problem}\label{equalm}
Determine $d(l,m;k) := d(m,m,\dots,m;k)$ for $l$ points in
$\mathbb P(V)$.
\end{problem}

Even in this somewhat reduced form, no general answer is known. In
fact, using Nagata's algorithm, one is naturally forced to
consider different $m_i$s even if one starts with equal
multiplicities, as we will see later. The case of equal
multiplicities has another nice feature: In Nagata's
counterexample one needs conditions of the form $\mathfrak a_m^k
\not= \mathfrak a_{km}$. Denote by $I_m = \bigoplus_k I_{m,k}$ the
homogeneous ideal of polynomials vanishing at our $l$ points to
equal multiplicity $m$. For $n = 2$, to show that $I_{m}^k \not=
I_{km}$ for certain arbitrarily large values of $m,k$, it
suffices to know that
\begin{equation}\label{nagatasineq}
  d(l,m;k)\not=0 \Rightarrow k > \sqrt{l}m
\end{equation}
Notice thats $d(l,m;k) \geq \binom{k+2}{2} -
l\binom{m+1}{2}$ is always true. This is just saying that the rank of the
equations imposed by the vanishing conditions is at most 
$l\binom{m+1}{2}$, which follows for example from our discussion
of problem~\ref{independence}. Using this one can show that
(\ref{nagatasineq}) cannot hold for all $k,m$, if we have $I_m^k =
I_{km}$ for $m$ sufficiently large. This leads us to another
\begin{problem}
For $n,l,\mathbf m$ fixed, what is the minimal degree $k_{\mathbf
m}$ such that $d(\mathbf m;k) > 0$?
\end{problem}
Nagata conjectured his condition~(\ref{nagatasineq}) to be true
for all $l \geq 10$, more precisely:
\begin{conjecture}[Nagata]Suppose $l \geq 10$ and $p_1,p_2,\dots,p_l$ are $l$
generic points in $\mathbb P^2$. If $f \in I_{m}$ then
$$\deg f > \sqrt{l}m$$
\end{conjecture}
The reason for requiring at least $10$ points is that the
assertion in Nagata's conjecture is false for all $l \leq 9$. One
might ask, why Steinberg was still able to use Nagata's method in
constructing a counterexample (\cite{STEIN}). The reason is,
that for nine points one still has $\deg f \geq \sqrt{9}m$ with
equality if and only if $f$ is a power of the unique cubic
determined by $9$ generic points in $\mathbb P^2$. Modifications
of Nagata's arguments then allow to conclude again that $I_{m}^k
\not = I_{km}$.

On the affirmative side, however, we have (cf.\ \cite{NAG59})
\begin{theorem}[Nagata, 1959]With the notation above, if $l = s^2 > 9$
is a square, then $f \in I_{m}$ implies $\deg f > sm$.
\end{theorem}

Before we go any deeper into the subject, we want to outline the
situation for $n=1$, which is actually the only dimension where
everything is clear. Suppose $p_1,\dots,p_l$ are generic points in
$V$. We may assume that $p_i = (1,c_i)$ for some $c_i \in F$,
hence all $p_i$ lie on the embedded copy $C = \{1\}\times F$ of $F$. If $f \in I_{\mathbf m,k}$ then $\phi =
f\left |_C\right .$ is an element of $F[t] = \mathcal O(F)$. Clearly $\phi(t) = g(t)\prod_i(t-c_i)^{m_i}$ with $\deg g
\leq k - \sum_i m_i$. Thus the dimension of $I_{\mathbf m;k} \left
| _C\right .$ is exactly $k+1 - \sum_i m_i$, if this is positive,
or zero otherwise. On the other hand, if $f\left |_C \right .$ is             
zero, then $f$ is zero. It follows that $d(\mathbf m;k) = k + 1 -
\sum_i m_i$ if this is positive, or zero otherwise, which is the
minimal possible. In other words, problem~\ref{independence} has a
positive answer in the case $n = 1$.

\section{The `dual' picture} \label{dual}

We now return to the original problem of determining
$d(m_1,m_2,\dots,m_l;k)$. Let $S = S(V) = \bigoplus_{k\geq 0}
S^d(V) = \bigoplus_{k\geq 0} S^d$ be the symmetric algebra of $V$,
which we think of as $\mathcal O(V)^* = \Hom_F(\mathcal O(V),F)$.
By our remark~\ref{vanishingorder}, we may replace the points
$p_1,\dots,p_l$ by elements of $V$ itself. So, for $p \in V$ and a
nonnegative integers $m$ and $k$, let $I_{p,m,k} = I_{Fp,m,k}$.
Then we have

\begin{lemma} For all $p \in V$ and all $m,k \geq 0$,
\begin{align}
I_{m,p,k} & = (S^{m-1}(V)p^{k-m+1})^\bot  \\ & = \{ f \in \mathcal
O(V)_k \mid \langle f, \alpha p^{k-m+1}\rangle = 0, \forall \alpha
\in S^{m-1}\}
\end{align}
\end{lemma}
Here $\langle \cdot, \cdot \rangle$ denotes the natural pairing
between  $\mathcal O(V)$ and $S(V)$.

\begin{proof} For any $v \in V$, let $D_v$ be the associated
derivation of $\mathcal O(V)$. Note that $\langle D_v(f), \lambda
\rangle = \langle f, v\lambda \rangle $ for all $v \in V, f \in
\mathcal O(V)_k, \lambda \in S^{k-1}$. To see this, without loss
of generality we may assume that $v = e_1$ the first coordinate
vector, then if $\mathcal O(V) = F[x_1,x_2,\dots,x_{n+1}]$ we have
$D_v = \frac{\partial}{\partial x_1}$. Since $S^{k-1}$ is spanned
by $(k-1)$th powers of elements of $v$, we are reduced to the case
$\lambda = w^{k-1}$ for some $w \in V$. Then $\langle
D_{e_1}(f),w^{k-1}\rangle = D_{e_1}(f)(w)$. Writing $f = f_0 +
f_1x_1 + f_2x_1^2 +\dots + f_kx^k$ with $f_i \in
F[x_2,\dots,x_{n+1}]$, we get $$\langle f,e_1w^{k-1}\rangle =
\langle f_0,e_1w^{k-1}\rangle +  \langle f_1x_1,e_1w^{k-1}\rangle
+ \dots + \langle f_kx_1^{k},e_1w^{k-1}\rangle,$$ which is
$0+f_1(w)+2f_2(w)x_1(w) +\dots +kf_k(w)x_1^{k-1}(w) =
D_{e_1}(f)(w)$.

But now we are done: By induction on $m$ we may assume that, $f
\in I_{p,m,k}$ if and only if $D_v(f)$ is orthogonal to
$S^{m-2}p^{k-m+1}$ for all $v$ and $f(p)  =\langle f,p^{k}\rangle
= 0$. In other words, if and only if $f$ is orthogonal to
$vS^{m-2}p^{k-m+1}$ for all $v$, which is our assertion.
\end{proof}

\begin{corollary}
For $l$ points we have
\begin{equation}
I_{\mathbf m,k} =  \bigl ( \sum_{i=1}^l S^{m_i-1}p_i^{k-m_i+1}
\bigr )^\bot
\end{equation}
\end{corollary}
As before $I_{\mathbf m,k}$ is the degree $k$-part of $I_{\mathbf
m} = I_{\mathbf p, \mathbf m}$.  The point is,
that if we put 
\begin{equation}
h(\mathbf p; \mathbf m; k) = \dim_F \bigl (
S^k/(\sum_iS^{k-r_i}p_i^{r_i})\bigr),
\end{equation}
then we have
\begin{equation}
d(\mathbf m; k) = h(\mathbf p; \mathbf m; k)
\end{equation}
with $r_i = k-m_i +1$, and the points $p_i$ sufficiently general. In other words, knowing $d(\mathbf m;k)$
for all $k$ and $\mathbf m$ is the same as knowing
\begin{equation}
h(\mathbf r;k) := \min\{ h(\mathbf p; \mathbf
r; k)| \mathbf p \in \mathbb P(V)^{l}\}
\end{equation}
for all $k$ and $\mathbf r = (r_1,r_2,\dots,r_l)$. 

\begin{remark}
It is a priori not clear whether generic points always
exist. But for any finite collection of  values for $\mathbf r$ it is
not hard to see that a minimizing choice of $\mathbf p$ exists such that
 $h(\mathbf p,\mathbf r;k) = h(\mathbf r; k)$ for all
$k$. 
\end{remark}

For fixed
$\mathbf r$ let $H_{\mathbf r}(q)$ be the formal power series
given by $ \sum_k h(\mathbf r;k)q^k$. Then $H_{\mathbf r}$ is the
Hilbert series of $S_{\mathbf r} := S/(Sp_1^{r_1}+Sp_2^{r_2}+\dots
+Sp_l^{r_l})$ for sufficiently generic points $p_i$. And if all $r_i$ equal $r$ ($1 \leq i \leq l$) we
abbreviate $H_{\mathbf r}$ by $H_{l,r}$ and $S_{\mathbf r}$ by
$S_{l,r}$.

\begin{definition} Let $f = \sum_i f_i q^i \in \mathbb Q\bigl[[q]\bigr]$ be a formal power
series. Then we set $$\trunc (f) = \sum_i \trunc (f_i)q^i$$ where
for $r \in \mathbb Q$
$$\trunc(r) = \begin{cases} r & \text{if} \quad r>0 \\ 0 & \text{otherwise}\end{cases}$$
\end{definition}

For $\mathbf r = (r_1,\dots,r_l)$ as above we set $$C_{\mathbf
r}(q) = \sum_k c(\mathbf r;k)q^k$$ where $$c(\mathbf r;k) =
\trunc\left ( \binom{k+n}{n}-\sum_{i=1}{l}\binom{k-r_i+2}{2}
\right )$$ is the virtual dimension of $S_{\mathbf r}^k$. Notice
that
\begin{equation}
C_{\mathbf r} =
\trunc\left(\frac{1-\sum_iq^{r_i}}{(1-q)^n}\right).
\end{equation}
If all the $r_i$ equal a given $r$ ($1 \leq i \leq l)$ we set
$C_{l,r} = C_{\mathbf r}$. Based on large scale computations for
values of $r$ between $10$ and $20$ it seems reasonable to
\begin{conjecture}
If $n = 2$ and $l\geq 10$, then for all $r$, $H_{l,r} = C_{l,r}$.
\end{conjecture}
For $l = 4$ and $l=9$ the assertion is true (despite the fact that
it is stated only for $l \geq 10$). This is a consequence of
Nagata's algorithm (cf.\ section~\ref{ninepts}). For $l=1,2$ it
fails for obvious reason: In these cases $S_r$ is not zero
dimensional as a ring, hence cannot have a polynomial as Hilbert
series. For $l = 3$ $S_r$ is the tensorproduct of three copies of $\mathbb
C[x]/(x^r)$, and so has $\frac{(1-q^r)^3}{(1-q)^3}$ as Hilbert series. For
$l=5,6,7,8$ the situation is more delicate. The conjecture fails for $l=5$ in
degree $4$ for $r=3$, for $l=6$ in degree $24$ for $r=15$. For $l=7$,
$h(l,r;k) \not = c(l,r;k)$ for $k = 42$ and $r=27$. Finally, $l=8$ fails in
degree $96$ with $r=63$.

An interesting test case for the conjecture is a situation where
$\binom{k+2}{2} = l\binom{k-r+2}{2}$. For example this occurs for
$l = 10$ for $k = 174$, $r = 120$. It does occur twice before, and the answer
there is $0 = h(10,r;k)$, which is well known. Using Groebner basis
techniques on a parallel computer system, we have
shown $h(10,r;174) = 0$, giving further evidence for the truth of the
assertion. The computation took several days, and so far it
seems to be the last number where once could actually compute the result. We
have added the resulting Hilbert series in the appendix below.

Returning to the general question in dimension $2$, we saw above,
that if $n=1$, then $H_{\mathbf r} = C_{\mathbf r}$. From this, it
is easy to see the following
\begin{lemma} For $n = 2$, $h(l,\mathbf r;k) = c(l,\mathbf r;k)$,
whenever $k > \sum_i m_i$.
\end{lemma}

\begin{proof} It suffices to show that the lemma is true for some points
$p_1,p_2,\dots,p_l \in \mathbb P^2$. We may choose the $p_i$
sitting on a line $L \subset \mathbb P^2$, such that they are
generic when considered as points of $L \cong \mathbb P^1$.
Suppose $L$ is given by an equation $l \in (F^2)^*$. In $S$ we may
therefore assume that $x,y,z$ are variables (on $(F^3)^*$)and that
$p_i \in F[x,y]_1 \subset S^1$ for all $i$. We look at the natural
injective restriction $r: S^k \rightarrow k[x,y]$, given by $f
\mapsto f(x,y,1)$, and with image all polynomials of degree less
or equal to $k$. Under this map, the $p_i$ remain homogeneous. It
follows that
 $$ \sum_i S^{k-r_i}p_i^{r_i} \cong \sum_i \bigoplus_{j = 0}^{k-r_i}r(p_i)^{r_i}F[x,y]_{j}$$
The latter equals
\begin{equation}\label{bigsum}
 \bigoplus_{k^\prime = 0}^k \sum_i r(p_i)^{r_i}F[x,y]_{k^\prime
 - r_i}
\end{equation}
Now $k > \sum m_i$ implies $k^\prime > \sum (k^{\prime} - r_i +
1)$. Since we know, that the the inner sum in \eqref{bigsum} is
direct under these circumstances, we are done.
\end{proof}

We conclude this section by proving the conjecture for $l=4$. So
suppose $p_1,p_2,p_3,p_4$ are four generic points of $V = F^3$. Because $\dim \mathcal O(V)_2 = 6$, any five points lie on a
quadric. So our four points lie, say, on $Q \subset \mathbb P(V)$,
defined by $f_0$ with $f_0$ homogeneous of degree $2$. Moreover
for generic points $f_0$ may be assumed generic also, so up to the
action of $GL_3$, all points may be chosen to have the form $p_i =
(1,t_i,t_i^2)$ and $f_0 = xz-y^2$. We have to show, that for fixed
$r>0$, we have $h(r;k) = c(r;k)$ for all $k$. We prove this by
induction on $r+ j = k$. Thus, we are investigating
$d(4,j+1;j+r)$. If $k = 0$, then either $r$ or $j$ is zero. For $r
= 0$ the result is trivial, and for $j=0,r=1$ the question is,
whether $4 = \dim \sum Fp_i^r$, which is true for the $p_i$
generic. Now suppose $h(r,r+j) = c(r;r+j)$ for all $r,j$ with $0
\leq r+j < k$, and we prove the result for $k$. As in the case
$n=1$ let $C $ be the embedded copy of $F$,
parameterized by $(1,t,t^2), t\in F$. For any $f \in I_{j+1,r+j}$,
let $\phi(t) = \phi_f(t) = f(1,t,t^2) \in F[t]$ be the restriction
of $f$ to $C$. Then it follows that
\begin{equation}
\phi(t) = g(t)\prod_{i=1}^4(t-t_i)^{j+1}
\end{equation}
with $\deg g \leq 2(r+j)-4(j+1)$. If $2(r+j)-4(j+1) <0$, there is
no such $g$ and $f\left |_C \right . = 0$, so $f = uf_0$ with $u
\in I_{j,r+j-2}$. If $r^\prime = r-1$ and $j^\prime = j-1$, it
follows that $2(r^\prime+j^\prime)-4(j^\prime+1) $ still is
negative, so by the very same argument $\phi_u(t) = 0$, and so on.
It follows $f = cf_0^d$ for some $d < j+1$ with $c$ linear or
constant. But $c$ vanishes on $C$, because $d <j+1$ and $j+1$ is
the required multiplicity for $f$, hence $f = 0$. On the other
hand it is easy to see, that $c(r,r+j) = 0$.

So suppose $2(r+j)\geq 4(j+1)$. In this case we have $$\dim
I_{j+1,r+j} \left | _C\right. \leq 2(r+j)-4(j+1) +1.$$ Using our
induction hypothesis this means that $h(4,r;r+j) \leq
c(4,r-1;r+j-2) + 2(r+j) - 4(j+1) +1$, because the space of
functions in $I_{j+1,r+j}$ restricting to zero on $C$ is just
$f_0P{j,r+j-2}$. But $c(4,r-1;r+j-2)$ is now nonnegative. Indeed
$2(r+j)-4(j+1) \geq 0$, so $r+j \geq 2j+2$, and
\begin{align*}
\binom{r+j}{2} & = \frac{(r+j)(r+j-1)}{2} \geq
\frac{(2j+2)(2j+1)}{2} \\
& = 4\frac{(j+1)(j+\frac{1}{2})}{2}  > 4\frac{(j+1)j}{2} =
\binom{j+1}{2}.
\end{align*}
It follows that $c(4,r-1;r+j-2) = \trunc\left(\binom{r+j-2 + 2
}{2}-4\binom{j-1 + 2}{2}\right)$ is actually positive. Hence,
\begin{align}
h(r;r+j) & \leq 2(r+j) - 4(j+1) + 1 +
\binom{r+j}{2}-4\binom{j+1}{2} \\
& = \binom{r+j+2}{2}-4\binom{j+2}{2} = c(r;r+j).
\end{align}
But certainly $h(4,r;r+j) \geq c(4,r;r+j)$ and the claim follows.

It should be remarked, that exactly the same argument works also
for nine points, except that one has to use a cubic instead of
our quadric $xz-y^2$ here. We will outline a proof in section~\ref{ninepts}. And the argument definitely does not work for
$5,6,7$ points.

\section{$n+2$ points in $\mathbb P^n$ and representations of $SL_2$}\label{onemore}

Now we turn to a different perspective. Suppose we are given
$p_1,p_2\dots,p_{l}$ generic points in $\mathbb P(V) =\mathbb
P^{n+1}$, and assume that $l \geq n+1$. It is clear that
$p_1,\dots,p_{n+1}$ my be taken to be linearly independent. In
this case, if we take $e_i \in p_i$ for $i = 1,2,\dots,n+1$ as a
basis for $V$, the elements
$e_1^{r_1},e_2^{r_2},\dots,e_{n+1}^{r_{n+1}}$ are a system
of parameters, and also a regular sequence of the ring $S(V)$. Set
\begin{equation}
M = S(V) / \left(\sum_{i=1}^{n+1}S(V)e_i^{r_i}\right)
\end{equation}
and view $M$ as an $S(V) \cong F[e_1]\otimes F[e_2] \otimes \dots
\otimes F[e_{n+1}] $-module, we get
\begin{equation}
M = F[e_1]/e_1^{r_1}F[e_1]\otimes F[e_2]/e_2^{r_2}F[e_2] \otimes
\dots \otimes F[e_{n+1}]/e_{n+1}^{r_{n+1}}F[e_{n+1}]
\end{equation}
with $e_i$ acting on the $i$th factor by multiplication and fixing
the other factors identically. On $F[e_i]/e_i^{r_i}F[e_i]$ we
choose a basis $v_{i,j}$ with $v_{i,j} = e_i^{r_i-j-1}$, $j =
0,2,\dots,r_i-1$. Then $e_i$ acts relative to this bases by
\begin{equation}
\begin{bmatrix}
0 & 1 & 0 & \cdots & 0 \\
0 & 0 & 1 & \cdots & 0 \\
\vdots & \vdots & \ddots & \ddots & \vdots \\
0 & 0 & \cdots & 0 & 1 \\
0 & 0 & \cdots & 0 & 0
\end{bmatrix}
\end{equation}
i.e.\ one nilpotent Jordan block. Define an endomorphism $h_i$ by
$h_i(v_{i,j}) = (r_i - 1 - 2j)v_{i,j}$, then there is exactly one
nilpotent operator $f_i$ with $[e_i,f_i] = h_i$ and $[h_i,f_i] =
-2f_i$. In other words $e_i,f_i,h_i$ form an $\mathfrak
sl_2$-triple and give rise to an $SL_2$-action on $M$, by
extending the action on the $i$th factor by the trivial
representation to all others. It follows $M$ is a  representation
of the product of $n+1$ copies of $SL_2$. Moreover
\begin{equation}
M = L^{r_1} \otimes L^{r_2} \otimes \dots \otimes L^{r_{n+1}}
\end{equation}
where for each nonnegative integer $r$, $L^r$ is the irreducible
$SL_2$-module of dimension $r$. To avoid clumsy notation we define
$L^0  = 0$, which corresponds to the case that one of $r_i = 0$
(and so $M = 0$).

For the rest of this section, we will assume that $l = n+2$, and
so we have one additional point $p = p_{n+2} \in \mathbb P(V)$.
Since we are interested in dimensions only it is clear, that we
may replace $p_1,p_2,\dots$ by $gp_1,gp_2,\dots$ where $g \in
GL_{n+1}$ is any element. Furthermore , since the $p_i$ are
generic, we may assume that $p = [a_1e_1 + a_2e_2 +\dots +
a_{n+1}e_{n+1}] \in \mathbb P(V)$ with all the $a_i$ nonzero. The
common stabilizer of our first $n+1$ points is the usual diagonal
torus of $GL_{n+1}$, so we may replace $p$ by $[e_1 + e_2 +\dots
+e_{n+1}]$, and set $e = e_1 + e_2 + \dots + e_{n+1}$. This last
definition makes sense also in $Lie((SL_2)^{n+1})$, so put $f =
f_1 + f_2 + \dots +f_{n+1}$ and $h = h_1 + h_2 + \dots +
h_{n+1}$. Therefor we get another $\mathfrak {sl}_2$-triple
corresponding to the diagonal $G$ in $(SL_2)^{n+1}$. Moreover the
action of $e \in Lie(G)$ on $M$ is the same as the action of the
\emph{point} $e \in V$.

For any vectorspace $Z$, on which $h$ acts, and any integer
$\lambda$, we define $Z[\lambda]$ to be the $\lambda$-eigenspace
of $h$ on $Z$. Moreover for any degree $k$, $M_k$ is the weight
space $M[\lambda]$ for $h$ for a suitable integer $\lambda$. Thus,
finding $H_{\mathbf r}$ is equivalent to finding the decomposition
of $M/e^rM$ as an $h$-module, where we put $r = r_{n+2}$ for
short.
\begin{lemma} Suppose $\mathbf r = (r_1,r_2,\dots,r_{n+2})$, then
\begin{equation}\label{hilbert-rep}
q^{-(r_1 + r_2 +\dots + r_{n+1}-{n+1})}H_{\mathbf r}(q^2) =
\sum_{\lambda}q^\lambda \dim (M/e^{r_{n+2}}M)[\lambda]
\end{equation}
where $\mathfrak {sl}_2$ acts via the tensor product action on $M =
L^{r_1}\otimes L^{r_2}\otimes \dots L^{r_{n+1}}$.
\end{lemma}
\begin{proof}
To keep things short, set $\lambda(k) = 2k - r_1 - r_2 \dots -
r_{n+1}+(n+1)$. Then $M[\lambda(k)] = M_k$ in our discussion of
$M$ as an $S(V)$-module above, and of course also
$(M/e^{r_{n+2}}M)[\lambda(k)] = (M/e^{r_{n+2}}M)_k$. It follows
that
\begin{equation}
q^{-(r_1 + r_2 +\dots + r_{n+1}-{n+1})}H_{\mathbf r}(q^2) =
\sum_{k}h(\mathbf r;k)q^{\lambda(k)},
\end{equation}
and the latter is also the righthand side of (\ref{hilbert-rep}).
\end{proof}
For $n = 1$ this gives a geometric interpretation of the
Clebsch-Gordan formula and for $n=2$ it gives an interpretation of the
$6j$-symbol. We also note that if $\mathbf r =
(r_1,r_2,\dots,r_{n+1},1)$, then finding $H_{\mathbf r}$ is
equivalent to finding the decomposition of $M$ into irreducibles.
Once this is done, it is a simple matter to replace the $1$ in the
$(n+2)$th position by $r_{n+2}$.

If all $r_i$ are equal, say, to $r$, then it is also possible to
deduce the multiplicities when $H_r$ is known. For example, what
we proved above for four points amounts to:
\begin{lemma}
$$L^r \otimes L^r \otimes L^r = \bigoplus_{j=0}^{r-1}
(j+1)L^{3(r-1)-2j+1} \oplus \bigoplus_{j=1}^{\left [ \frac{r-1}{2}
\right ]} (r-2j)L^{r-2j}$$
\end{lemma}

In the rest of this section we write $n$, the dimension of
$\mathbb P(V)$, as a subscript to avoid any confusion. Note that
$H_{n,(r_1,r_2,\dots,r_{n+1},1)} = H_{n,(r_1,\dots,1,r_{n+1})}$
and also $S/e_{n+1}S = F[e_1,\dots,e_n]$. Thus we get
$$H_{n,(r_1,\dots,r_{n+1},1)} = H_{n-1,(r_1,\dots,r_{n+1})}.$$
Here on the right hand side everything takes place in one
dimension less. In the special case of $n = 2$, it follows, that
for the coefficients of $H_{2,(r_1,r_2,r_3,1)}$, we have
\begin{equation}
\label{oneless}h_{2}((r_1,r_2,r_3,1);k) = \trunc \bigl(k+1 -
\sum_{k-r_i \geq 0}(k - r_i + 1) \bigr )
\end{equation}
On the other hand, the left hand side of (\ref{oneless}) is $$\dim
\bigl ( (L^{r_1} \otimes L^{r_2} \otimes L^{r_3}) / e( L^{r_1}
\otimes L^{r_2} \otimes L^{r_3})[-r_1 - r_2 - r_3 + 3 + 2k]\bigr
).$$ Summarizing, this implies the following generalization to the
Clebsch-Gordan formula.
\begin{proposition}
Let $r_1,r_2,r_3,l$ be positive integers. If $l \equiv r_1 + r_2 +
r_3 \mod 2$, then with $$k = \frac{r_1 + r_2 + r_3 - 3 - l}{2},$$
we have
$$\dim \Hom_{SL_2}(L^{l+1}, L^{r_1} \otimes L^{r_2} \otimes L^{r_3}) = \trunc \Bigl ( k+ 1 -  \sum_{k - r_i \geq 0} (d-r_i + 1)\Bigr ).$$
Otherwise $\Hom_{SL_2}(L^{l+1}, L^{r_1} \otimes L^{r_2} \otimes
L^{r_3})$ is $(0)$.
\end{proposition}

In the last section of this paper, we will give an algorithm, that computes
$H_{\mathbf r}$ for at most nine points in $\mathbb P^2$. By what was said
above, to compute the multiplicities in $L^{r_1}\otimes \dots \otimes
L^{r_4}$, one has to know $H_{3,(r_1,r_2,\dots,r_4,1)} =
H_{2,(r_1,r_2,\dots,r_4)}$. We may assume that $r_1 \leq r_2 \leq r_3 \leq r_4$, and we  borrow the following result from section~\ref{ninepts}:

\begin{proposition}\label{borrow} Suppose $2k + 3 \leq  r_1 + r_2 + r_3$ and $0
\leq r_1 \leq r_2 \leq r_3 \leq r_4 \leq k$, then $h_2(r_1,r_2,r_3,r_4;k) =
c_2(r_1,r_2,r_3,r_4;k) > 0$.
\end{proposition}

As we will see now, this is almost everything one has to know for decomposing
a fourfold tensorproduct. We are immediately reduced to the case where $2k+ 3
> r_1 + r_2 + r_3$ or $r_4 > k$. In all other cases, proposition~\ref{borrow} gives the correct answer. 
As before set $\lambda(k) = -(r_1 + r_2 + r_3) + 3 + 2k$, the corresponding weight of $h$ in degree $k$.
If $r_4 \leq \lambda(k)$, then the transformation rule in the next section asserts that $h_2(r_1,r_2,r_3,r_4;k) = 0$. On the other hand, if $r_4 > k$, the answer is given by $h_2(r_1,r_2,r_3;k)$, the coefficient of $q^{k}$ in
$$\frac{(1-q^{r_1})(1-q^{r_2})(1-q^{r_3})}{(1-q)^3}.$$
Finally, if $\lambda(k) \leq r_4\leq k$, the dimension of $(e^{r_4}M)[\lambda(k)]$ is the sum of multiplicities of $L^p$ in $M$, where $p$ ranges over the set $p-1 \equiv \lambda(k) \mod 2$, $p > 2r_4 - \lambda(k)$. If we write $p = \lambda(k) + 1 + 2j$, this multiplicity is $h_1(r_1,r_2,r_3;k + j + 1)$.
Thus, in this case \[h_2(r_1,r_2,r_3,r_4;k) = h_2(r_1,r_2,r_3;k) - \sum_{j = r_4 - \lambda(k)}^\infty h_1(r_1,r_2,r_3;k+j+1).\]

\section{The transformation rule} \label{transform}
We now turn to the central point of our discussion. Our aim is to
redevelop Nagata's algorithm from our point of view. There is a
universal rule which simply says that $h(\mathbf r;k) = h(\mathbf
r^\prime, k^\prime)$ for certain $\mathbf r^\prime$, $k^\prime$.
The algorithm depends on the fact that for nine or less points it
is possible to see exactly when $k^\prime < k$, and then to handle
the case $k^\prime \geq k$ directly. We will come to that in a
moment.

First we adapt our approach from the last section to more than
$n+2$ points. Suppose we are given $l$ generic points in $\mathbb
P(V)$. As before we may identify the first $n+1$ points with the
coordinate lines and choose representatives $e_1, e_2,\dots,
e_{n+1}$. Moreover, we may assume that all the other $p_i$ have
nonzero coordinates, i.e. $p_i = [a_{i1}e_1 + a_{i2}e_2 +\dots +
a_{in+1}e_{n+1}]$, with all the $a_{ij}$ nonzero. As before we
interpret $e_i$ as a certain element of the $i$th factor in
$(\mathfrak {sl}_2)^{n+1}$. Also $e,f,h$ are defined as in the
last section. Given this suppose $t$ is an element of the maximal
torus $T \subset (SL_2)^{n+1}$ consisting of diagonal matrices.
Then $t$ corresponds to an $n+1$-tuple $(t_1,t_2,\dots,t_{n+1})$
of points in $F^*$, and we have
$$\Ad(t)e = t_1^2e_1 + t_2^2e_2 + \dots + t_{n+1}^2e_{n+1}.$$
Taking residue classes in $$M = S(V) /e_1^{r_1}S(V) +
e_2^{r_2}S(V) + \dots + e_{n+1}^{r_{n+1}}S(V)$$ we see that $e_j =
\sum_i a_{ji}e_i$ for $j > n+1$ is of the form $\Ad(t_j)e$ for a
suitable $t_j \in T$. Thus,
\begin{equation}
e_j^{r_j}M = t_je^{r_j}M \qquad j = n+2,n+3,\dots,l.
\end{equation}
Of course, $t_je^{r_j}M$ and also $e^{r_j}M$ are $h$-stable, and
even isomorphic as $h$-modules.
Notice also, that $e_i,f_i,h_i = h$ form an $\mathfrak{sl}_2$-triple where $f_i =
\Ad(t_i)f$ (and $h = \Ad(t_i)h$).

\begin{theorem}\label{t-rule}
If $l > n+1$ then,
\begin{equation}\label{t-rule-eq}
h(\mathbf r; k ) = h(\mathbf r^\prime; k^\prime)
\end{equation}
with $k^\prime = r_1 + r_2 + \dots + r_{n+1} - n - 1 - k$,
$\mathbf r^\prime = (r_1,r_2,\dots,r_{n+1},
r_{n+2}^\prime,r_{n+3}^\prime,\dots,r_l^\prime)$, and $r_i^\prime
=  r_1 + r_2 + \dots + r_{n+1} - n - 1 + r_i - 2k$.
\end{theorem}
\begin{remark}
As mentioned before, in the case $n = 2$  the theorem is
equivalent to a result of Nagata {\rm(\cite{NAG60_II})}, which he
proves by methods of classical projective geometry.

On the other hand, the proof we present here uses only elementary
representation theory of $SL_2$.
\end{remark}

We will need the following Lemma in the proof.

\begin{lemma}\label{t-rule-help}
Let $L^p$ be an irreducible $\mathfrak{sl}_2$-representation with highest
weight $p-1$ and let $e,f,h$ be a standard triple inside $\mathfrak{sl}_2$. 
Then for all integers $r,\lambda$ with $r \geq 0$ we have:
If $r > \lambda$, then
\begin{equation}\label{e_is_f_eq}
(e^rL^p)[\lambda] = (f^{r-\lambda}L^p)[\lambda].
\end{equation}
If on the other hand $r \leq \lambda$, then
\begin{equation}\label{e_is_all_eq}
L^p[\lambda] = (e^r L^p)[\lambda].
\end{equation}
Here for an $h$-module $V$, $V[\lambda]$ is the $h$-weight space of
weight $\lambda$. 
\end{lemma}

\begin{proof}
First we deal with the case $r >\lambda$.
Notice that in (\ref{e_is_f_eq}) both sides of the equation are either zero or
one-dimensional and in the latter case they equal $L^p[\lambda]$. It therefore
is enough to show that the dimensions agree, i.e.\ that they are nonzero for
the same values of $r$ and $\lambda$. If $\lambda$ and $p-1$ do not have the
same parity, then $L^p[\lambda]$ is zero and therefore (\ref{e_is_f_eq}) is
obviously true. Hence it is safe to assume that $\lambda$ and $p-1$ both are
simultaneously odd or even.

In this case the left hand side of~\ref{e_is_f_eq} is nonzero if and only if
\begin{equation}\label{e-eq}
-p + 1 +2r \leq \lambda \leq p-1.
\end{equation}
Similarly the right hand side is nonzero if and only if 
\begin{equation}
\label{f-eq}
-p+1 \leq \lambda \leq p-1-2(r-\lambda).
\end{equation}
Of course both these equations are equivalent: assuming (\ref{e-eq}), it
follows immediately that $-p+1 \leq \lambda$ because $r \geq 0$. And
$p-1-2(r-\lambda) = p-1-2r+2\lambda \geq -\lambda + \lambda = \lambda$ due to
the left hand side of (\ref{e-eq}), hence (\ref{f-eq}) holds.
A completely analogous argument shows the other implication, that is, if
(\ref{f-eq}) holds, (\ref{e-eq}) is true as well.

Finally assume that $r \leq \lambda$. In this case, if $\lambda$ is a weight
of $L^p$, so is $\lambda - 2r$. But then $L^p[\lambda] = e^r(L[\lambda - 2r])
= (e^rL^p)[\lambda]$, hence the claim.
\end{proof}

Using this we are able to prove the theorem.

\begin{proof}[Proof of theorem~\ref{t-rule}]
We are interested in $$\dim\, \Bigl(M/\sum_{i>n+1}
e_i^{r_i}M\Bigr)[\lambda(k)]$$ where $\lambda(k) = -r_1-r_2-\dots
-r_{n+2}+n+1 +2k$ is the $h$-weight associated to degree $k$.
First we will consider the trivial cases. Suppose $k^\prime$ is
negative. Thus, the right hand side of (\ref{t-rule-eq}) is zero
(by convention). On the other hand, this means that $k$ is greater
than the highest weight of $M$, and so $\lambda(k)$ is greater
than this weight, too. Thus $M[\lambda(k)] = (0)$, and
(\ref{t-rule-eq}) holds.

We may therefore assume, that $k^\prime \geq 0$. Notice that $r_i^ \prime =
r_i - \lambda(k)$. If there is an index $i> n+1$
such that $r_i^\prime \leq 0$ then $\lambda(k) \geq r_i$. Consider the action
of the $i$th $\mathfrak{sl}_2$-triple on $M$. For each irreducible submodule
$L^p \subset M$ we are in the second case of Lemma~\ref{t-rule-help} which
then asserts that $L^p[\lambda(k)]= (e_i^{r_i}L^p)[\lambda(k)]$. Thus
$M[\lambda(k)] = e_i^{r_i}M[\lambda(k)]$ and the left hand side of
(\ref{t-rule-eq}) is zero. But the right hand side is zero as well because if
$r_i^{\prime} = 0$ this is obvious and if $r_i^{\prime} < 0$ this is
convention.

It remains to treat the case that all $r_i^{\prime}$ are strictly positive.
For each $i > n+1$ we are then in the first case of Lemma~\ref{t-rule-help}. That is, looking at the $i$th $\mathfrak{sl}_2$-triple acting on
$M$, then for each irreducible $L^p\subset M$ we have 
$$(e_i^{r_i}L^p)[\lambda(k)] = (f_i^{r_i^\prime}L^p)[\lambda(k)].$$
Of course, this applies to all of $M$ as well and we conclude, that for
each $i$
$$(e_i^{r_i}M)[\lambda(k)] = (f_i^{r_i^\prime}M)[\lambda(k)].$$

Moreover,
$$\lambda(k) = r_1 + r_2 + \dots + r_{n+1} - (n+1) - 2k^\prime = -\lambda(k^\prime).$$

For any $\mathfrak{sl}_2$-module $L$ there is an involution $\theta$
satisfying $\theta(ex) = f\theta(x)$, and $\theta(hx) = -h\theta(x)$. Thus
there is an involution $\Theta$ of $M$, satisfying $\Theta(e_ix) =
f_i\Theta(x) $ and $\Theta(h_ix) = -h_i\Theta(x)$ for $1 \leq i \leq n+1$
($\Theta$ is obtained by tensoring the $\theta$s of the individual factors
$L^{r_i}$.) Clearly this implies 
$$\Theta \Bigl(\sum_{i>n+1} (f_i^{r_i^ \prime}M)[\lambda(k)]\Bigr) = \sum_{i>n+1}(e_i^{r_i^\prime}M)[\lambda(k^\prime)]$$ because $\lambda(k^\prime) = -\lambda(k)$,
and the theorem now follows.
\end{proof}
It is quite remarkable that one is able to deduce this
transformation rule using only elementary representation theory of
$SL_2$. In contrast to that, Nagata's method for $n = 2$ was to
show that under a quadratic transformations $T$ of $\mathbb P^2$,
the global sections of a certain linear system
$\mathcal L$ are in one to one correspondence with the sections of $T\mathcal L$.

\section{$9$ points in $\mathbb P^2$ and Nagata's algorithm} \label{ninepts}
In this section we will develop our version of Nagata's algorithm
to determine the numbers $h(\mathbf r;k)$ for less than nine
points in $\mathbb P^2$.

So let $p_1,\dots,p_l$ be generic points of $\mathbb P^2$ with $3
\leq l \leq 9$. First we state the algorithm and then prove that
it terminates after finitely many steps. If we consider points in
another dimension than $2$, we will indicate that with a
subscript, e.g.\ we will write $h_1(\mathbf r;k)$ for $l$ generic
points in $\mathbb P^1$.
\begin{algorithm}\label{alg} With the notation above, if the test at stage
$i$ fails we go to $i+1$, otherwise we are done or start at $1$
again.
\begin{enumerate}
\item If $l = 3$, $h(r_1,r_2,r_3;k)$ is the coefficient of $q^k$
in $$ \prod_{i=1}^{3}\frac{1-q^{r_i}}{1-q}.$$
\item Put $r_1,r_2,\dots,r_l$ in increasing order. If $r_1 \leq
0$, $h(\mathbf r;k) = 0$.
\item If $r_1 = 1$, then $$h(\mathbf r;k) = h_1(r_2,r_3,\dots,r_l; k) =
\trunc\bigl(k + 1 - \sum_{i=2}^l \trunc(k - r_i + 1)\bigr).$$
\item If $r_l > k$, then $h(r_1,\dots,r_l;k) =
h(r_1,r_2,\dots,r_{l-1};k)$.
\item If $2k + 3 > r_1 + r_2 + r_3$, we set $r_i^\prime =  r_1 +
r_2 + r_3 - 3 + r_i - 2k$ and $k^\prime = r_1 + r_2 + r_3 -3 - k$,
and $$h(\mathbf r;k) =
h(r_1,r_2,r_3,r_4^\prime,\dots,r_{l}^\prime;k ^\prime).$$
\item $h(\mathbf r;k) = h(r_1-2,r_2-2,\dots,r_l - 2;k-3) + 3k -
\sum_i (k-r_i+1)$.
\end{enumerate}
\end{algorithm}
We assert that this is indeed an algorithm. Suppose, $1$ applies.
Then it is well known that the Hilbert series of
$F[x_1,x_2,x_3]/(x_1^{r_1},x_2^{r_2},x_3^{r_3})$ has the asserted
form, since it is the tensor product of the $F[x]/(x^{r_i})$, as
noted above. Hence, we may assume that the $r_i$ are ordered, and
we are in step $2$. If $r_1$ is less than zero, we already
observed that by convention $h( \mathbf r;k)$ is zero. If $r_1 =
0$, then $\sum S^jp_j^{r_j} \subset S^k$, and it also follows that
$h(\mathbf r;k) = 0$. Hence, we are in step $3$. But if $r_1 = 1$,
the claimed equation is obvious. And in step $4$, if $r_l > k$,
$S^lp_l^{r_l}$ does not contribute to degree $k$, so $r_l$ may be
dropped. We are now at step $5$. Here, $l \geq 4$, and if $2k + 3>
r_1 + r_2 +r_3$ then $k^\prime < k$ and $r_i^\prime < r_i$. Hence
step $5$ reduces the degree and the $r_i$ for $i > 3$, and because
of the transformation rule, it preserves the value of $h$.

It thus remains to look at step $6$. For this we may assume that
$l \geq 4$, $2k +3 \leq r_1 + r_2 + r_3$, and $2 \leq r_1 \leq r_2
\leq r_3 \leq \dots \leq r_l \leq k$. Note that we must have $r_1
\geq 3$, otherwise $2k + 3 \leq r_1 + r_2 + r_3$ implies $r_3 >
k$. If $r_1 = 3$, then all of the $r_i$ equal $k$ for $i > 1$, by
the same arguments. If we can prove that in these circumstances we
have $$ h(\mathbf r;k ) = c(\mathbf r;k ) = \binom{k+2}{k} -
\sum_{i=1}^l\binom{k-r_i + 2}{2},$$ then step $6$ will yield a
correct answer. For this we will need the following lemma.

\begin{lemma} Assume that $2k + 3 \leq r_1 + r_2 + r_3$, and $r_1
\leq r_2  \leq \dots \leq r_l \leq k$. Then $c(\mathbf r;k) > 0$.
\end{lemma}
\begin{proof}Without loss of generality, $l = 9$. Put $m_i = k - r_i + 1$. Then $m_1 \geq m_2 \geq
\dots \geq m_l \geq 0$. Moreover, $$m_1 + m_2 + m_3 = 3k + 3 - r_1
- r_2 - r_3 \leq k.$$ We have to show that
\begin{equation} \label{c_ge_0}
\binom{k+2}{2} \geq \sum_i \binom{m_i+1}{2} + 1
\end{equation}
The right hand side of (\ref{c_ge_0}) is 
\begin{align*} \sum_i
\frac{m_i(m_i+1)}{2}+ 1 & = \frac{1}{2}\sum_i (m_i^2 + m_i) + 1\\ & \leq 
\frac{m_1^2 + m_2^2 + 7m_3^2 + m_1 + m_2 + 7m_3 + 2}{2}.
\end{align*}
The left hand side is at least
\begin{align*}
&\frac{(m_1+m_2+m_3 + 1)(m_1+m_2+m_3+2)}{2}\\
&= \frac{m_1^2+m_2^2 +m_3^2 + 2m_1m_2 + 2m_1m_3 + 2m_2m_3 + 3m_1 +
3m_2+3m_3+2}{2}\\
&\geq \frac{m_1^2+m_2^2+7m_3^2+m_1+m_2+7m_3+2}{2}.
\end{align*}
The last inequality is because $m_3 \leq m_2,m_1$. The lemma
follows.
\end{proof}
Returning to step $6$ of the algorithm, we will prove a
slightly stronger statement: 

\begin{proposition}
Suppose $0 \leq r_1 \leq r_2 \leq
\dots \leq r_l \leq k$ with $2k + 3 \leq r_1 + r_2 + r_3$, then
$h(\mathbf r;k) = c(\mathbf r;k) > 0$.
\end{proposition}
 Obviously,
proposition~\ref{borrow} in section~\ref{onemore} is a special case. The will be by induction on
$k$.

\begin{proof}
Restricting our attention to points of the form $(1,t_i,t_i^3)$
with $t_1 + t_2 + \dots + t_l \not = 0$ and using arguments
similar to those in section~\ref{dual} for four points, we see
with $\mathbf r_0 = (r_1 - 2,r_2-2,\dots,r_l-2)$ and $k_0 = k - 3$
that
\begin{equation}\label{upper_bound}
h(\mathbf p; \mathbf r;k) \leq 3k - \sum_{i = 1}^l (k - r_i + 1) +
h(\mathbf p;\mathbf r_0;k_0).
\end{equation}
Note that we are using special points here.
But these
are as generic as any others, as we will see in a moment. One might expect the right hand side of (\ref{upper_bound})
to be too small by one. However, any homogeneous polynomial of
degree $k$ has, when restricted to the curve $(1,t,t^2)$, a zero
coefficient in front of $t^{3k-1}$. Since our points satisfy
$\sum_i t_i \not = 0$ it is easy to see that this is really a non
trivial condition for polynomials vanishing to the given orders.
Thus, the possible dimension of the restricted space is one less.

It is clear that $2k_0 + 3 \leq (r_1 - 2) + (r_2 - 2) + (r_3 - 2
)$, hence we are in a similar situation as when we started. There are three cases to consider: First, assume that $r_i < k$
for all $i$. Then also $r_i - 2 \leq k_0$, and induction yields
the result: $c(\mathbf r; k) = 3k - \sum_i (k - r_i + 1) + c(\mathbf
r_0;k_0) > 0$, which one sees easily by straightforward computations.
Thus, we may consider the second case: there is $i > 3$
with $r_i = k$. Since the $r_i$s are ordered, this holds starting
from a $j$ throughout to $l$. When computing $h(\mathbf r_0;k_0)$
the points $p_j,\dots,p_l$ thus may be dropped. Again, induction
and the last lemma assert that $h(\mathbf p;\mathbf r_0;k) = c(\mathbf
r_0^*;k_0)
> 0$ where $\mathbf r_0^* = (r_1-2,r_2-2,\dots,r_{j-1}-2)$. On the other
hand, obviously $0 < c(\mathbf r;k) = c(\mathbf r^*;k) -
\sum_{i=j}^l 1$. But the latter is just the right hand side of
(\ref{upper_bound}) together with the induction hypothesis. Again,
$\mathbf r^*$ is gotten by dropping $r_j,\dots,r_l$.

Finally, we have to consider the case when $j \leq 3$. We are now
reduced to prove the following: $$3k - (k-r_1+1) - (k-r_2+1) - 1 +
h(\mathbf p; r_1-2,r_2-2,k-2;k-3) - (l-3) \leq c(\mathbf r;k)$$ This
follows easily from the fact, that two points always generate
independent conditions, if $r_1 + r_2 \geq k+3 \geq k$:
$$(x_1^{r_1}F[x_1,x_2,x_3])_k \cap (x_2^{r_2}F[x_1,x_2,x_3])_k = (0)$$
if $r_1 + r_2 > k$. Thus, $$h(\mathbf p; r_1-2,r_2-2,k-2;k-3) =
h(\mathbf p; r_1-2,r_2-2;k-3)$$ and the latter equals
$$\binom{k-1}{2}-\binom{k-r_1+1}{2}-\binom{k-r_2+1}{2},$$ because
$r_1+r_2 \geq k+3$, hence $r_1-2 + r_2-2 \geq k-3$, and the claim
follows. In particular, our $l$ points are generic, and we may
drop $\mathbf p$ in $h$.
\end{proof}

We can now prove:
\begin{theorem}
Suppose $l = 4$ or $9$. Then $H_r = C_r$. In other words $h(r;k) =
c(r;k)$ for all $r$ and $k$.
\end{theorem}
\begin{proof}
We have already seen the case $l = 4$, so we stick with the case
$l = 9$. All multiplicities are equal now. If we are able to prove
that $2k + 3 > 3r$ implies that $h(r;k) = 0$, we are done, by what
we have said above, since we know that the conjectural formula is
true if $2k + 3 \leq 3r$ and $r \leq k$.

So suppose $2k + 3 > 3r$. This means also that $3k \leq 9(k-r+1)$.
If we look at our curve of the form $(1,t,t^3)$ and assume our
points sitting on it, it follows that all restrictions vanish.
Thus, $h(r,k) = h(r-2,k-3)$. But also $2(k-3)+3 > 3(r-2)$, and we
may conclude that $h(r-2,k-3) = 0$, provided $r < k$. If $r = k$,
then $2k + 3 > 3k$ implies $k=0,1,2$, and so $h(k;k) = 0$.
\end{proof}
\begin{remark}
The theorem is also a consequence of the results of Steinberg in \cite{STEIN}.
\end{remark}

\appendix

\newpage

\section*{Appendix: The Hilbert series $H_{10,120}$ of $S_r$ for $10$
  sufficiently generic points in $\mathbb P^2$ with $r = 120$}

$$375 q^{173} + 741 q^{172} + 1098 q^{171} + 1446 q^{170} + 1785
  q^{169} + 2115 q^{168} + 2436 q^{167} +$$ $$
 2748 q^{166} + 3051 q^{165} + 3345 q^{164} + 3630 q^{163} + 3906 q^{162} +
 4173 q^{161} + 4431 q^{160} +$$ $$ 4680 q^{159} + 4920 q^{158} + 5151 q^{157}
 + 5373 q^{156} + 5586 q^{155} + 5790 q^{154} + 5985 q^{153} +$$ $$ 6171
 q^{152} + 6348 q^{151} + 6516 q^{150} + 6675 q^{149} + 6825 q^{148} + 6966
 q^{147} + 7098 q^{146} +$$ $$ 7221 q^{145} + 7335 q^{144} + 7440 q^{143} +
 7536 q^{142} + 7623 q^{141} + 7701 q^{140} + 7770 q^{139} +$$ $$ 7830 q^{138}
 + 7881 q^{137} + 7923 q^{136} + 7956 q^{135} + 7980 q^{134} + 7995 q^{133}
 +8001 q^{132} +$$ $$ 7998 q^{131} + 7986 q^{130} + 7965 q^{129} + 7935
 q^{128} + 7896 q^{127} + 7848 q^{126} + 7791 q^{125} + $$$$ 7725 q^{124} +
 7650 q^{123} + 7566 q^{122} + 7473 q^{121} + 7371 q^{120} + 7260 q^{119} + 7140 q^{118} + $$$$7021 q^{117} + 6903 q^{116} + 6786 q^{115} + 6670 q^{114} +
6555 q^{113} + 6441 q^{112} +  6328 q^{111} + $$$$ 6216 q^{110} + 6105 q^{109} + 5995 q^{108} +
 5886 q^{107} + 5778 q^{106} + 5671 q^{105} +5565 q^{104} + $$$$ 5460 q^{103}
 + 5356 q^{102} + 5253 q^{101} + 5151 q^{100} + 5050 q^{99} + 4950 q^{98} +4851 q^{97} + $$$$4753 q^{96} + 4656 q^{95} + 4560 q^{94} + 4465 q^{93} + 4371
 q^{92} +4278 q^{91} + 4186 q^{90} + 4095 q^{89} +$$$$ 4005 q^{88} + 3916 q^{87} + 3828 q^{86} + 3741 q^{85} + 3655 q^{84} + 3570 q^{83} + 3486 q^{82} + 3403 q^{81} +$$$$ 3321 q^{80} + 3240 q^{79} + 3160 q^{78} + 3081 q^{77} + 3003 q^{76} + 2926 q^{75} + 2850 q^{74} + 2775 q^{73} +$$$$ 2701 q^{72} + 2628 q^{71} + 2556 q^{70} + 2485 q^{69} + 2415 q^{68} + 2346 q^{67} + 2278 q^{66} + 2211 q^{65} +$$$$ 2145 q^{64} +2080 q^{63} + 2016 q^{62} + 1953 q^{61} + 1891 q^{60} + 1830 q^{59} + 1770 q^{58} + 1711 q^{57} +$$ $$ 1653 q^{56} + 1596 q^{55} + 1540 q^{54} + 1485 q^{53} + 1431 q^{52} + 1378 q^{51} + 1326 q^{50} + 1275 q^{49} +$$$$ 1225 q^{48} + 1176 q^{47} + 1128 q^{46} + 1081 q^{45} + 1035 q^{44} + 990 q^{43} +946 q^{42} + 903 q^{41} +$$$$ 861 q^{40} + 820 q^{39} + 780 q^{38} + 741 q^
{37} + 703 q^{36} + 666 q^{35} + 630 q^{34} + 595 q^{33} + 561 q^{32} +$$$$ 528 q^{31} + 496 q^{30} + 465 q^{29} + 435 q^{28} + 406 q^{27} +
 378 q^{26} + 351 q^{25} + 325 q^{24} + 300 q^{23} +$$$$ 276 q^{22} + 253 q^{21}
 + 231 q^{20} + 210 q^{19} + 190 q^{18} + 171 q^{17} + 153 q^{16} + 136 q^{15}
 + 120 q^{14} + $$$$105 q^{13} + 91 q^{12} + 78 q^{11} + 66 q^{10} + 55 q^{9} +
 45 q^{8} + 36 q^{7} + 28 q^{6} + 21 q^{5} + 15 q^{4} + 10 q^{3} + $$$$ 6 q^{2} +
 3 q^{1} + 1$$


\begin{thebibliography}{9}
\bibitem{A} J. Alexander, \textit{Singunlarit\'es imposables en position
    g\'en\'erale \`a une hypersurface projective} Compositio Math. \textbf{68}
    (1988), no. 3, 305 - 354
\bibitem{AH92} J. Alexander, A. Hirschowitz, \textit{La m\'ethode d'Horace
    \'eclat\'ee: application \`a l'interpolation en degr\'e quatre},
    Invent. Math. \textbf{107} (1992), 585 - 602
\bibitem{CM2000} C. Ciliberto, R. Miranda, \textit{Linear systems of plane
    curves with base points of equal multiplicity},
    Trans. Amer. Math. Soc. \textbf{352} (2000), no. 9, 4037 - 4050
\bibitem{DF99} D. Daigle, G. Freudenburg,  \textit{A counterexample to
    Hilbert's fourteenth problem in dimension $5$}, J. Algebra \textbf{221}
    (1999), no. 2, 528-535
\bibitem{NAG59} M. Nagata, \textit{On the $14$-th problem of Hilbert},
    American J. Math. \textbf{81} (1959), 766 - 772
\bibitem{NAG60_I} M. Nagata, \textit{On rational surfaces. I. Irreducible
    curves of arithmetic genus $0$ or $1$}, Mem. Coll. Sci. Univ. Kyoto Ser. A
  Math. \textbf{32} (1960), 351 - 370
\bibitem{NAG60_II} M. Nagata, \textit{On rational surfaces. II},
  Mem. Coll. Sci. Univ. Kyoto Ser. A \textbf{33} (1960/1961), 271 - 293
\bibitem{STEIN} R. Steinberg, \textit{Nagata's example}, Algebraic groups and
  Lie groups, 375 - 384, Austral. Math. Soc. Lect. Ser., 9, Camebridge
  Univ. Press, Camebridge, 1997


\end{thebibliography}
\end{document}